\documentclass{amsart}
\usepackage{enumerate}
\usepackage[dvips]{graphicx}
\usepackage{color}
\usepackage{pstricks}
\usepackage{pst-node}
\newtheorem{theo}{Theorem}[section]

\newtheorem{lem}[theo]{Lemma}

\newtheorem{prop}[theo]{Proposition}

\newtheorem{defn}[theo]{Definition}

\theoremstyle{definition}

\def\a{{\a}}
\def\a{{\mathfrak a}}

\def\C{\mathbb C}

\def\DD{\mathcal D}
\def\CC{\mathcal C}

\def\N{\mathbb N}

\def\R{\mathbb R}

\def\Z{\mathbb{Z}}
\begin{document} 

\title[Generalized LUE and interacting particles with a wall]{Generalized Laguerre Unitary Ensembles\\ and \\ an interacting particles model with a wall}

\author{Manon Defosseux}
\address{Laboratoire de Math\'ematiques Appliqu\'ees \`a Paris 5, Universit\'e Paris 5, 45 rue des  Saints P\`eres, 75270 Paris Cedex 06.}
\email{manon.defosseux@parisdescartes.fr}

\begin{abstract} We introduce and study a new interacting particles model with a wall and two kinds of interactions -  blocking and pushing - which maintain particles in a certain order.  We show that it involves a random matrix model. 
\vspace{0.3cm}\\
{\sc R\'esum\'e.}  Nous introduisons et \'etudions un nouveau mod\`ele de particules en interaction avec mur. Les particules de ce mod\`ele se bloquent et se poussent les une les autres, afin de rester dans un ordre donn\'e. Nous montrons que ce mod\`ele est li\'e \`a un mod\`ele de matrices alŽatoires. 

\end{abstract}
\maketitle
\section{Interacting particles model}\label{model}
Let us consider $k$ ordered particles evolving in discrete time on the positive real line with interactions that maintain their orderings. The particles are labeled in increasing order from $1$ to $k$. Thus  for $t\in \N$,  we have $$0\le X_1(t)\le \dots\le X_k(t),$$ where  $X_i(t)$ is the position of the $i^{\textrm{th}}$ particle at time $t\ge 0$.  Particles are initially all at $0$. The particles jump at times $n-\frac{1}{2}$ and $n$, $n\in \N^*$. Let us consider two independent families  $$(\xi(i,n-\frac{1}{2}))_{i=1,\dots,k;\, n\ge 1}, \quad \textrm{and } \quad (\xi(i,n))_{i=1,\dots,k;\, n\ge 1},$$   of independent random variables  having  an exponential law of mean $1$.  For convenience, we suppose that there is a static particle which always stays at $0$. We call it the $0^{\textrm{th}}$ particle,  and denote $X_0(t)$ its position at time $t\ge 0$.

\noindent \underline{At time $n-1/2$}, for $i=1,\dots, k$, in that order, the $i^{th}$ particle  tries to jump to the left according to a jump of size $\xi(i,n-\frac{1}{2})$ being blocked by the old position of the $(i-1)^{th}$ particle.    In other words :
\begin{itemize}
\item If $X_1(n-1)-\xi(1,n-\frac{1}{2})<0$, then the $1^{st}$ particle  is blocked by $0$, else it jumps at $X_1(n-1)-\xi(1,n-\frac{1}{2})$, i.e. $$X_1(n-\frac{1}{2})=\max(0,X_1(n-1)-\xi(1,n-\frac{1}{2})).$$  
\item For $i=2,\dots,k$, if $X_i(n-1)-\xi(i,n-\frac{1}{2})<X_{i-1}(n-1)$, then the $i^{th}$ particle  is blocked by $X_{i-1}(n-1)$, else it jumps at $X_i(n-1)-\xi(i,n-\frac{1}{2})$, i.e. $$X_i(n-\frac{1}{2})=\max(X_{i-1}(n-1),X_i(n-1)-\xi(i,n-\frac{1}{2})).$$
\end{itemize}
\underline{At time $ n$}, particles jump successively to the right according to an exponentially distributed jump of mean $1$, pushing all the particles to their right. The order in which the particles jump is given by their labels.  Thus  
for $i=1,\dots,k$,  if $X_{i-1}(n)>X_{i}(n-\frac{1}{2})$ then the $i^{th}$ particle is pushed before jumping, else it jumps from $X_{i}(n-\frac{1}{2})$ to  $X_{i}(n-\frac{1}{2})+\xi(i,n)$:    
  $$X_{i}(n)=\max(X_{i-1}(n),X_{i}(n-\frac{1}{2}))+\xi(i,n).$$ 
  
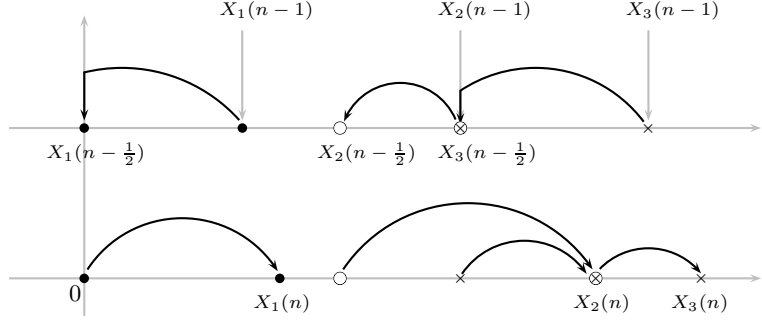
\begin{figure}[h!]
\begin{pspicture}(-9.5,5)(6,0.4)
\psline[linecolor=lightgray]{->}(-8,1)(2,1) \psline[linecolor=lightgray]{->}(-8,3)(2,3)\psline[linecolor=lightgray]{->}(-7,0.5)(-7,4.5)

\psline[linecolor=lightgray]{<-}(0.5,3.1)(0.5,4.3)\put(0.2,4.5){\scriptsize{$X_{3}(n-1)$}}
\psline[linecolor=lightgray]{<-}(-2,3.1)(-2,4.3)\put(-2.3,4.5){\scriptsize{$X_{2}(n-1)$}}
\psline[linecolor=lightgray]{<-}(-4.9,3.1)(-4.9,4.3)\put(-5.2,4.5){\scriptsize{$X_{1}(n-1)$}}

\psarc{-}(-1,2){1.8}{37}{124}\psline{->}(-2,3.5)(-2,3.05)
\psarc{->}(-2.8,2.8){0.8}{23}{159}

\psdots[dotstyle=o,dotsize=5pt](-3.6,3)(-2,3)
\psdots[dotstyle=x](0.5,3)(-2,3)

\psarc(-6.5,1.7){2.1}{42}{104}\psline{<-}(-7,3.1)(-7,3.75) 
 \psdots(-4.9,3)(-7,3)

\put(-2.3,2.6){\scriptsize{$X_{3}(n-\frac{1}{2})$}}\put(-3.9,2.6){\scriptsize{$X_{2}(n-\frac{1}{2})$}}\put(-7.5,2.6){\scriptsize{$X_{1}(n-\frac{1}{2})$}}

\psdots[dotstyle=o,dotsize=5pt](-3.6,1)(-0.2,1)
\psdots[dotstyle=x](-2,1)(-0.2,1)(1.2,1)
 \psdots(-7,1)(-4.4,1)

\psarc{<-}(-5.7,0.3){1.5}{32}{147} 
\psarc{<-}(-1.9,0){2}{34}{146} 
\psarc{<-}(-1.15,0.5){1}{34}{146} 
\psarc{<-}(0.5,0.5){0.9}{39}{136} 

\put(-7.2,0.7){\small{$0$}}\put(0.8,0.6){\scriptsize{$X_{3}(n)$}}\put(-0.5,0.6){\scriptsize{$X_{2}(n)$}}\put(-4.75,0.6){\scriptsize{$X_{1}(n)$}}

\end{pspicture}
\caption{An exemple of blocking and pushing interactions between times $n-1$ and $n$ for $k=3$.}
\label{model}
\end{figure}

There is another  description  of the dynamic of this model which is equivalent  to the previous.  At each time $n\in \N^*$, each particle successively attempts to jump first to the left then to the right, according to independent exponentially distributed  jumps of mean 1. The order in which the particles jump is given by their labels.  At time $n\in \N^*$, the $1^{st}$  particle  jumps to the left being blocked by $0$ then immediately  to the right, pushing the $i^{\textrm{th}}$ particles, $i=2,\dots,k$. Then the second particle jumps to the left, being blocked by $\max(X_1(n-1),X_1(n))$, then to the right, pushing the $i^{\textrm{th}}$ particles, $i=3,\dots,k$, and so on. In other words, for $n\in \N^*$, $i=1,\dots,k$,
\begin{align*}
X_i(n)=\max\big(X_{i-1}(n),X_{i-1}(n-1),X_{i}(n-1)-\xi^{-}(i,n)\big)+\xi^+(i,n),
\end{align*}
where   $(\xi^-(i,n))_{i=1,\dots,k;\, n\in \N}$, and $(\xi^+(i,n))_{i=1,\dots,k;\, n\in \N}$ are two independent families of independent random variables having  an exponential law of mean $1$. 

\section{Results}
 Let us denote ${\mathcal M}_{k,m}(\R)$ the  real vector space of $k \times m$ real  matrices.  We put on it the Euclidean structure defined by the scalar product, 
$$\langle M, N \rangle = \mbox{tr}(MN^*), \quad M,N\in {{\mathcal M}_{k,m}(\R)}.$$
Our choice of the Euclidean structure above defines a notion of standard Gaussian variable on  ${{\mathcal M}_{n,m}(\R)}$.  We write ${\mathcal A}_{k}$ for the set $\{M\in \mathcal{M}_{k,k}(\R): M+M^*=0\}$ of antisymmetric $k\times k$ real matrices, and  $i{\mathcal A}_{k}$  for the set $\{iM: M\in \mathcal{A}_k\}$. Since a matrix in $i\mathcal{A}_k$ is Hermitian, it has real eigenvalues $\lambda_1 \geq \lambda_2 \geq \cdots \geq \lambda_k$. Morever, antisymmetry implies that  $\lambda_{k-i+1}=-\lambda_{i}$, for $i=1,\cdots, [k/2]+1$, in particular $\lambda_{[k/2]+1}=0$ when $k$ is odd.
 
Our main result is that the positions of the particles of our interacting particles model can be interpreted as  eigenvalues of a random walk on the set of matrices $i\mathcal{A}_{k+1}$. 

\begin{theo}\label{theorem} Let $k$ be a positive integer and $(M(n),n\ge 0)$, be a discrete process  on $i\mathcal{A}_{k+1}$ defined by $$M(n)=\sum_{l=1}^{n}Y_l \begin{pmatrix} 0 &  i  \\
-i &  0       
\end{pmatrix}Y_l^*,$$ where the $Y_l$'s are independent standard Gaussian variables in $\mathcal{M}_{k+1,2}(\R)$. For $n\in \N$, let  $\Lambda_1(n)$ be the largest eigenvalue of $M(n)$. Then the processes $$(\Lambda_1(n),n\ge 0)\textrm{  and } (X_{k}(n),n\ge 0),$$ have the same distribution.
\end{theo}
For a matrix $M\in i\mathcal{A}_{k+1}$ and $m\in\{1,\dots,k+1\}$, the main minor of order $m$ of $M$ is the submatrix $$\{M_{ij}, \, 1\le i,j\le m\} .$$
\begin{theo}\label{minor} Let $M(n),n\ge 0$, be a discrete process  on $i\mathcal{A}_{k+1}$, defined as in theorem \ref{theorem}.  For $n\in \N$, $m=2,\dots,k+1$, we denote $\Lambda_1^{(m)}(n)$, the largest eigenvalue of the main minor of order $m$ of $M(n)$. Then for each fixed $n\in \N^*$, the random vectors
$$(\Lambda_1^{(2)}(n),\dots,\Lambda_1^{(k+1)}(n)) \quad \textrm{ and } \quad (X_{1}(n),\dots,X_{k}(n)),$$
have the same distribution.
\end{theo}
Let us notice that there already exists a version of theorem \ref{theorem}  with particles jumping by one (see section 2.3 of \cite{WarrenWindridge}) and a continuous version which involves reflected Brownian motions with a wall and Brownian motions conditioned to never collide with each other or the wall (see \cite{WarrenetCie}).

There is a variant of our model with no wall and no left-jumps which has been extensively studied (see for instance Johansson \cite{Johansson}, Dieker and Warren \cite{DiekerWarren}, or Warren and Windridge \cite{WarrenWindridge}). It involves random matrices from the Laguerre Unitary Ensemble. Indeed, in that case, the position of the rightmost particle has the same law as the largest eigenvalue of the process $(M(n),n\ge 0)$ defined by $$M(n)=\sum_{l=1}^nZ_lZ_l^*, \quad n\ge 0,$$ where the $\sqrt{2}\,Z_i$'s are independent standard Gaussian variables on  $\C^k$. Let us mention that this model without a wall is equivalent to a maybe more famous one called the TASEP (totally asymmetric simple exclusion process). In this last model, we consider  infinitely many ordered particles evolving on $\Z$ as follows. Initially there is one and only one particle on each point of $\Z_-$. Particles  are equipped with independent Poisson clocks of intensity 1 and each of them jumps by one to the right only if its clock rings and the point just to its right is empty. Particles are labeled by $\N^*$ from  the right to the left. Then, the time needed for the $k^{th}$ particles to make $n$ jumps  is exactly the position of the rightmost particle at time $n$ in the model without a wall and exponential right jumps.

\section{consequences}
Thanks to the previous results we  can deduce some properties of the interacting particles from the known properties concerning the matrices $M(n),n\ge 0$. For instance, the next proposition follows immediately  from theorem \ref{theorem} and theorem 5.2 of \cite{Defosseux}.
\begin{prop}  For $n\ge [\frac{k+1}{2}]$, the distribution function of $X_{k}(n)$ is given by

$$\mathbb{P}(X_{k}(n)\le t)=\det\big(\int_0^tx^{2j+i+n-[\frac{k+1}{2}]-3+1_{\{k\textrm{ is even}\}}}\, e^{-x}\, dx\big)_{[\frac{k+1}{2}]\times [\frac{k+1}{2}]},\quad t\in \R_+.$$
\\

\end{prop}

Moreover, as $n$ goes to infinity, the process $(\frac{1}{\sqrt{n}}M([nt]),t\ge 0)$ converges in distribution to a Brownian motion on $i\mathcal{A}_{k+1}$. Since our interacting particles model converges to a model of reflecting Brownian motions, proposition 2 of \cite{WarrenetCie} follows from theorem \ref{theorem}.

 \section{proofs}
Theorem \ref{theorem} is proved by proposition \ref{intertwining}. Theorem \ref{minor} follows from propositions \ref{corintertwining} and \ref{lawminor}.
Let us denote by $\phi$ the function from $\R$ to $\R$, defined by $\phi(x)=\frac{1}{2}exp(-\vert x\vert)$, $x\in \R$, and  consider the random walk $S_n,n\ge 1$, on $\R$, starting from $0$, whose increments have  a density equal to $\phi$. The next three lemmas are elementary. 
\begin{lem}  \label{absoluterw} The process $(\vert S_n\vert)_{n\ge 1}$ is a Markov process, with a transition density $q$ given  by 
$$q(x,y)=\phi(x+y)+\phi(x-y), \quad x,y\in \R_+.$$
\end{lem}
\proof This is a simple computation which holds for any symmetric random walk. $\Box$
\begin{lem}  \label{so2} Let us consider a   process $M(n),n\ge 0,$ on $i\mathcal{A}_2$, defined by $$M(n)=\sum_{l=1}^{n}Y_l \begin{pmatrix} 0 &  i  \\
-i &  0       
\end{pmatrix}Y_l^*, \quad n\ge 0,$$ where the $Y_l$'s are independent standard Gaussian variables in $\mathcal{M}_{2,2}(\R)$. Then the process of the only positive eigenvalue of $(M(n), n\ge 0)$ is Markovian with transition density given by $q$. 
\end{lem}
\proof The $Y_l$'s are $2\times 2$ independent real random matrices whose entries are   independent standard Gaussian random variables on $\R$. Let us write $Y_l=\begin{pmatrix} Y_{l,1} &  Y_{l,2}    \\
Y_{l,3} &  Y_{l,4}        
\end{pmatrix}$.  The matrix $Y_l \,\omega(1)Y_l^*$ is equal to $$\begin{pmatrix} 0 & i(Y_{l,4}Y_{l,1} -Y_{l,2}Y_{l,3} )   \\
 i(Y_{l,2}Y_{l,3} -Y_{l,1}Y_{l,4})  &  0       
\end{pmatrix}.$$ For $\alpha\in \R$, we have $$\mathbb{E}( e^{-i\alpha(Y_{l,4}Y_{l,1} -Y_{l,2}Y_{l,3} )})=\frac{1}{1+\alpha^2}.$$
Thus, the random variables $Y_{l,4}Y_{l,1} -Y_{l,2}Y_{l,3} $, $l=1,\dots,n$, are independent, with   a density equal to $\phi$. We conclude using lemma \ref{absoluterw}. $\Box$\\

\begin{lem} \label{rrandomwalk} Let $r$ be a real number. Let us consider $(\xi^-_n)_{n\ge 1}$ and $(\xi^+_n)_{n\ge 1}$ two independent families of independent random variables having an exponential law of mean $1$. The Markov process $Z(n), n\ge 1$, defined by $$Z(n)=\max(Z(n-1)-\xi^-_n,r)+ \xi^+_n, n\ge 1,$$ has transition density 
$$p_r(x,y)=\phi(x-y)+e^{2r}\phi(x+y), \quad x,y\ge r.$$ 
\end{lem}
\proof This is a simple computation. $\Box$\\

Let us notice that when $r=0$, the Markov process $(Z(n), n\ge 0)$, describes the evolution of the first particle. As its transition kernel $p_0$ is the same as the transition kernel $q$ defined in lemma \ref{absoluterw},    theorem \ref{theorem} follows when $k=1$ from lemma \ref{so2}. The general case is more complicated.

\begin{defn} We define the function $d_k:\R^{[\frac{k+1}{2}]}\to \R$ by 
\begin{itemize}
\item when $k=2p$, $p\in \N^*$,
$$
d_k(x)=c_k^{-1}\prod_{1\le i<j\le p}(x^2_i-x^2_j)\prod_{i=1}^px_i,$$
 \item when $k=2p-1$, $p\in \N^*$,
$$d_k(x)= \left\{
    \begin{array}{ll}
      1 & \mbox{if } p= 1 \\
          c_k^{-1}\prod_{1\le i<j\le p}(x^2_i-x^2_j) & \mbox{else,}
    \end{array}
\right.$$
where $$c_k=2^{[\frac{k}{2}]}\prod_{1\le i<j\le p}(j-i)(k+1-j-i)\prod_{  1\le i\le p}(p+\frac{1}{2}-i)^{1_{\{k=2p\}}}.$$ 
\end{itemize}

\end{defn}

The next proposition gives the transition density of the process of eigenvalues of the process $M(n),n\ge 0$.  For the computation, we need a generalized Cauchy-Binet identity (see for  instance Johansson  \cite{Johansson2}).  Let $(E,\mathcal{B},m)$ be a measure space, and let  $\phi_{i}$ and $\psi_{j}$, $1\le i,j\le n$, be measurable functions such that the $\phi_{i}\psi_{j}$'s are  integrable. Then  \begin{align} \label{CauchyBinet}
\det\Big(\int_{E}\phi_{i}(x)\psi_{j}(x)dm(x)\Big)=\frac{1}{n!}\int_{E^{n}}\det\big(\phi_{i}(x_{j}) \big)\det \big(\psi_{i}(x_{j}) \big)\prod_{k=1}^n dm(x_{k}). 
\end{align}
We also need  an identity which expresses interlacing conditions with the help of a determinant (see  Warren \cite{Warren}).  For $x, y\in \mathbb{R}^{n}$ we write $x \succeq y$ if $x$ and $y$ are interlaced, i.e.\@
$$x_{1}\ge y_1\ge x_2 \ge \cdots\ge x_n \ge y_n  $$ and we write $x \succ y$ when
 $$ x_{1}> y_1> x_2 > \cdots> x_n> y_n. $$
When $x\in \mathbb{R}^{n+1}$ and $y\in\mathbb{R}^n$ we add the relation $y_n \geq x_{n+1}$ (resp. $y_n> x_{n+1}$). Let $x$ and $y$ be two vectors in $\mathbb{R}^{n}$ such that $x_{1}> \cdots >x_{n}$ and $y_{1}> \cdots >y_{n}$. Then \begin{align}  1_{x\succ y}=\det(1_{\{x_{i}>y_{j}\}})_{n\times n}\label{interlacing}.\end{align}
\\

For $k\ge 1$, we denote by $\CC_k$ the subset of $\R^{[\frac{k+1}{2}]}$ defined by $$\CC_k=\{x\in\R^{[\frac{k+1}{2}]}: x_{1} > \dots >x_{[\frac{k+1}{2}]} >0\}.$$

\begin{prop}
Let us consider the   process  $(M(n),n\ge 0)$, defined as in theorem \ref{theorem}. For $n\in \N$, let $\Lambda(n)$ be the first $[\frac{k+1}{2}]$ largest eigenvalues of $M(n)$, ordered such that $$\Lambda_1(n)\ge \dots\ge \Lambda_{[\frac{k+1}{2}]}(n)\ge 0.$$ Then $(\Lambda(n),n\ge 0)$, is a Markov process with a transition density $P_k$ given by 
$$P_k(\lambda,\beta)= \frac{d_k(\beta)}{d_k(\lambda)}\det(\phi(\lambda_i-\beta_j)+(-1)^{k+1}\phi(\lambda_i+ \beta_j))_{1\le i,j\le [\frac{k+1}{2}]}, \quad \lambda,\beta\in \CC_k.$$

\end{prop}

\proof The Markov property follows from the fact that the matrices $$Y_l \begin{pmatrix} 0 &  i  \\
-i &  0       
\end{pmatrix}Y_l^*, \quad l\in \N^*,$$ are independent and have an invariant distribution for the action of the orthogonal group by conjugacy. Proposition 4.8 of \cite{Defosseux} ensures that the transition density with respect to the Lebesgue measure of the positive eigenvalues of $M(n),n\ge 1$, is given by 
\begin{itemize}

\item  when  $k=2p$,
 \begin{align*}
 P_k(\lambda,\beta)=\frac{d_k(\beta)}{d_k(\lambda)} I_p(\lambda,\beta),\quad \lambda,\beta\in \CC_k ,\end{align*}

\item  when  $k=2p-1$,
\begin{align*}
 P_k(\lambda,\beta)=\frac{d_k(\beta)}{d_k(\lambda)} \frac{1}{2}(e^{-\vert \lambda_p-\beta_{p}\vert }+e^{-( \lambda_p+\beta_{p}) }) I_{p-1}(\lambda,\beta), \quad \lambda,\beta\in \CC_k,\end{align*}
 \end{itemize}
 where $$I_p(\lambda,\beta)= \left\{
    \begin{array}{ll}
      1 & \mbox{if } p= 0 \\
          \int_{\mathbb{R}_+^{p}} 1_{\{\lambda,\beta\succ z\}}e^{-\sum_{i=1}^{p}(\lambda_i+\beta_i-2z_i)}\,dz& \mbox{else.}
    \end{array}
\right.$$
When $k$ is even, using identity (\ref{interlacing}), we write $1_{\{\lambda,\beta\succ z\} } e^{-\sum_{i=1}^{p}(\lambda_i+\beta_i-2z_i) }$ as $$\det(1_{z_i< \lambda_j}e^{-(\lambda_j-z_i)})_{p\times p}\det(1_{z_i<\beta_j}e^{-(\beta_j-z_i)})_{p\times p},$$
and use the Cauchy-Binet identity to get the proposition. 

\noindent When $k$ is odd, we introduce the measure $\mu$ on $\R$, defined by $\mu=\delta_0+m$, where $\delta_0$ is the Dirac measure at $0$ and $m$ is the Lebesgue measure on $\R$. We have the identity 
$$ \phi(x-y)+\phi(x+y)=\, \int_\R1_{[0,x\wedge y]}e^{-(x+y-2z)}\,d\mu(z).$$ Thus using the Cauchy-Binet identity with the measure $\mu$ we get that the determinant of the proposition is equal to 
$$\frac{1}{p!}\int_{\R^{p}}\det\big(1_{[0,\lambda_j]}(z_{i})e^{-(\lambda_j-z_i)} \big)\det \big(1_{[0,\beta_j]}(z_{i})e^{-(\beta_j-z_i)} \big)\prod_{m=1}^p\, d\mu(z_{m}).$$ Using identity (\ref{interlacing}), we obtain that it is equal to
$$\int_{\R^{p}}1_{\lambda,\beta\succ z}e^{-\sum_{i=1}^p(\lambda_i+\beta_i-2z_i)}\prod_{m=1}^p\, d\mu(z_{m}).$$ We integrate over $z_p$ in the last integral and use the fact that the coordinates of $\beta$ and $\lambda$ are strictly positive to get the proposition.  $\Box$\\

The next proposition gives the transition density of the Markov process $(X(n), n\ge 0)$ defined in section \ref{model}. For $m\ge 0$, we denote $\phi^{(m)}$ the $m^{\textrm{th}}$ derivative of $\phi$, and we define the function $\phi^{(-m)}$, by $$\phi^{(-m)}(x)=(-1)^m\int_{x}^{+\infty}\frac{1}{(m-1)!}(t-x)^{m-1}\phi(t)\,dt, \quad x\in \R.$$ We easily obtain that 
\begin{align*} 
\phi^{(m)}(x) &= \left\{
    \begin{array}{ll}
        \frac{1}{2}(-1)^me^{-x} & \mbox{if } x\ge 0 \\
          \frac{1}{2}e^{x} & \mbox{else}.
    \end{array}
\right.\\
\phi^{(-m)}(x) &= \left\{
    \begin{array}{ll}
       \phi^{(m)}(x) & \mbox{if } x\ge 0 \\
         -\sum_{i=1}^{[\frac{m+1}{2}]}x^{m-(2i-1)} +\phi^{(m)}(x)& \mbox{else}.
    \end{array}
\right.
\end{align*} 
\\

For $k\ge 2$, we denote by $\DD_k$ the subset of $\R^{k}$ defined by $$\DD_k=\{x\in\R^{k}: 0< x_1< x_2< \dots<x_{k}\}.$$
\begin{prop} \label{kernelparticles}The Markov process $X(n)=(X_1(n),\dots,X_k(n)), n\ge 0$, has a transition density  $Q_k$ given by
$$Q_k(y,y')=\det(a_{i,j}(y_i,y_j'))_{1\le i,j\le k}, \quad y,y'\in \bar{\DD}_k, $$
where for $x,x'\in \R$, $$a_{i,j}(x,x')=(-1)^{i-1}\phi^{(j-i)}(x+x')+(-1)^{i+j}\phi^{(j-i)}(x-x').$$
\end{prop}
\proof Let us show the proposition by induction on $k$.  For $k=1$, the equality holds by lemma  \ref{so2}. Suppose that it is true for $k-1$. We write $C_1,\dots, C_k$, and $L_1,\dots, L_k$ for the columns and the rows of the matrix of which we compute the determinant. There are two cases :
\begin{itemize}
\item  If $y'_{k-1}\ge y_{k-1}$ then for $i=1,\dots,k-1$, $y_i\le y_{k-1}'\le y_k'$.  A quick calculation shows that all the  components  of the  column $C_k+e^{y_{k-1}'-y'_k}C_{k-1}$ are equal to zero except the last one, which is equal to $p_{y'_{k-1}}(y'_{k-1}\vee y_k,y_k')$, where  $p_r$ is defined in lemma \ref{rrandomwalk} for $r\in \R$.
\item If $y'_{k-1}\le y_{k-1}$ then for $i=1,\dots,k-1$, $y_i'\le y_{k-1}\le y_k$. We replace the last line $L_k$ by the  line $L_k-e^{y_{k-1}-y_k}L_{k-1}$ having all its components equal to zero except the last one, which is equal to $p_{y_{k-1}}(y_k,y_k')$.  
\end{itemize}
Then we conclude developing the determinant according to its last column in the first case or to its last  row  in the second one, and using the induction property. $\Box$\\

We have now all the ingredients  needed to prove theorems \ref{theorem} and \ref{minor}. For $\lambda$ an element of the adherence $\Bar{\CC}_k$ of $\CC_k$, we denote GT$_k(\lambda)$  the subset of $\R^{[\frac{k+1}{2}][\frac{k+2}{2}]}$ defined by \begin{eqnarray*}\textrm{GT}_k(\lambda)=\{(x ^{(2)},\cdots, x ^{(k+1)}): x ^{(k+1)}=\lambda, x ^{(i)}\in \mathbb{R}_+^{[\frac{i}{2}]},\,   x ^{(i)}   \succeq    x^{(i-1)}  , 3\leq i\leq k+1\}.\end{eqnarray*} 
 We let $$\textrm{GT}_k=\cup_{\lambda\in \bar{\CC}_k}\textrm{GT}_k(\lambda).$$ If $(x ^{(2)},\cdots, x ^{(k+1)})$ is an element of GT$_k,$ then $(x_1^{(2)},\cdots, x_1^{(k+1)})$ belongs to the adherence $\bar{\DD}_{k}$ of $\DD_{k}$. Thus we define L$_k$ as the Markov kernel on $\bar{\CC}_k\times \bar{\DD}_{k}$ such that for $\lambda\in \bar{\CC}_k$, the probability measure L$_k(\lambda,.)$ is the  image   of the uniform probability measure on GT$_k(\lambda)$  by the projection  $p:\textrm{GT}_k\to \DD_{k}$ defined by $$p((x ^{(2)},\cdots, x ^{(k+1)}))=(x_1^{(2)},\cdots,x_1^{(k+1)}),$$ where $(x ^{(2)},\cdots, x ^{(k+1)})\in \textrm{GT}_k.$ For $\lambda\in \CC_k$, the volume of $\textrm{GT}_k(\lambda)$ is given by $d_k(\lambda)$. Thus  L$_k(\lambda,.)$ has a density with respect to the Lebesgue measure on GT$_k$ given by  $$\frac{1}{d_k(\lambda)}1_{x\in \textrm{GT}_k(\lambda)}, \quad x\in\R^{[\frac{k+1}{2}][\frac{k+2}{2}]}.$$
Rogers and Pitman proved in \cite{RogersPitman} that it is sufficient to show that the intertwining  (\ref{inter}) holds, to get the equality in law of the processes $(\Lambda_1(n),n\ge 0)$ and $(X_{k}(n),n\ge 0)$.  So theorem \ref{theorem} follows from proposition \ref{intertwining}.

\begin{prop}\label{intertwining}
\begin{align}\label{inter}
L_kQ_{k}=P_kL_k
\end{align}
\end{prop} 
\proof The proof is the same as the one of proposition 6 in  \cite{WarrenetCie}. We use the determinantal expressions for $Q_{k}$ and $P_k$ to show that both sides of equality (\ref{inter}) are equal to the same determinant. For  this we use that  the coefficients $a_{i,j}$'s given in proposition \ref{kernelparticles} satisfy for $x,x'\in \R_+$,
\begin{align*}
&a_{i,j}(x,x')=\int_x^{+\infty}a_{i-1,j}(u,x')\, du,\\
&a_{i,j}(x,x')=-\int_{x'}^{+\infty}a_{i,j+1}(u,x')\, du,\\
&a_{2i,2j}(x,0)=0,\quad a_{2i,2i-1}(0,x)=1,\quad a_{2i,j}(0,x)=0, \, 2i\le j
\end{align*}
The computation of the left hand side of (\ref{inter}) rests on the first identity.  The computation of the right hand side rests on the others. $\Box$\\

Let us notice that there are many reasons to think that we could construct an enlarged process $\big((Y^{(2)}(n),\dots,Y^{(k+1)}(n)), n\ge 0\big),$ living on GT$_k$, such that the processes $(Y^{(k+1)}(n),n\ge 0)$, and $\big((Y_1^{(2)}(n),\dots,Y_1^{(k+1)}(n)), n\ge 0\big),$ would have  respectively the same law as $(\Lambda(n),n\ge 0)$ and $(X_1(n),\dots,X_{k}(n),n\ge 0)$.   This  would imply theorem \ref{theorem}. 

The measure L$_k(0,.)$ is the Dirac measure at the null vector of GT$_k(0)$. Thus, the following proposition is an immediate consequence of proposition \ref{intertwining}.

\begin{prop}\label{corintertwining}
\begin{align*}
Q_{k}^n(0,.)=P_k^nL_k(0,.)
\end{align*}
\end{prop} 
Keeping the same notations as in theorem \ref{minor}, we have the following proposition, from which theorem \ref{minor} follows.
\begin{prop} \label{lawminor} $P_k^nL_k(0,.)$ is the law of the random variable $$(\Lambda_1^{(2)}(n),\dots,\Lambda_1^{(k+1)}(n))).$$
\end{prop} 
\proof The density of the positive eigenvalues $\Lambda(n)$ of $M(n)$ is given by $P_k^n(0,.)$.   Then the proposition follows immediately from theorem 3.4 of   \cite{Defosseux}. $\Box$

\section{concluding remarks}
As recalled, the model with a wall that we have introduced  is a variant of another one with no wall and no left-jumps. For this last model, the proofs of the analogue results as those of theorems \ref{theorem} and \ref{minor}  rest on some combinatorial properties of Young tableaux. Indeed, Young tableaux are used to describe the irreducible representations of the unitary group. The matrices from the Laguerre Unitary Ensemble belong to the set of Hermitian matrices which is, up to a multiplication by the complex $i$, the Lie algebra of the Unitary group. Their laws are  invariant for the action of the unitary group by conjugacy. It is a general result that  the law of their eigenvalues  can be deduced from some combinatorial properties of the irreducible representations of the unitary group.   

In our case, the distribution of the eigenvalues of $(M(n),n\ge 0)$ can be deduced from combinatorial properties of the irreducible representations of the orthogonal group   (see \cite{Defosseux} for details). Many  combinatorial approaches  have been developed to describe these representations. Among them we can mention the orthogonal tableaux and the analogue of the Robinson Schensted algorithm   for the  orthogonal group (see Sundaram \cite{Sundaram}), or  more recently those  based on the very general theory of crystal graphs of Kashiwara \cite{Kashiwara}. None of them seems to lead to the interacting particles model with a wall that we have introduced. It would be interesting to find what kind of tableau involved in the description of irreducible representations of the orthogonal group would lead to this model.    

\end{document}